\documentclass[12pt,a4paper]{article} 

\def\endofps{EndOfTheIncludedPostscriptMagicCookie}
\chardef\other=12
\newwrite\psdumphandle 
\outer\def\psdump#1{\par\medbreak
  \immediate\openout\psdumphandle=#1
  \copytoblankline}
\def\copytoblankline{\begingroup\setupcopy\copypsline}
\def\setupcopy{\def\do##1{\catcode`##1=\other}\dospecials
  \catcode`\\=\other \obeylines}
{\obeylines \gdef\copypsline#1
  {\def\next{#1}%
  \ifx\next\endofps\let\next=\endgroup %
  \else\immediate\write\psdumphandle{\next} \let\next=\copypsline\fi\next}}
\outer\def\closepsdump{
  \immediate\closeout\psdumphandle}

\psdump{hardy.eps}
statusdict begin/waittimeout 0 def end{/MGXPS_2.8_EPS begin}stopped{/MGXPS_2.8_EPS
200 dict def}if MGXPS_2.8_EPS begin/bd{bind def}bind def/xd{exch def}bd/ld{load
def}bd/M/moveto ld/m/rmoveto ld/L/lineto ld/l/rlineto ld/w/setlinewidth
ld/n/newpath ld/P/closepath ld/tr/translate ld/gs/gsave ld/gr/grestore
ld/lc/setlinecap ld/lj/setlinejoin ld/ml/setmiterlimit ld/TM/concat
ld/F/fill ld/bz/curveto ld/_m matrix def/_cm matrix def/_d 0 def/_o
0 def/_pb 0 def/_hb -1 def currentscreen/_P xd dup/_A xd/__A xd dup/_F
xd/__F xd/srgbc{_neg 0 ne{3{neg 1 add 3 1 roll}repeat}if setrgbcolor}bd/ssc{gs
0 0 M 0 _yp neg _xp 0 0 _yp 3 D closepath fill gr 0 0 0 srgbc}bd/Cne{/_bs
xd/_neg where{/_neg get/_neg xd pop}{pop 0/_neg xd}ifelse dup 0 ne{setflat}{pop}ifelse
dup -1 eq{pop}{pop pop}ifelse/_rz xd/_yp xd/_xp xd/_a 0 def 0 ne{/_a
90 def -90 rotate _xp _rz div 72 mul neg 0 tr}if 0 _yp _rz div 72 mul
tr 72 _rz div dup neg scale/#copies xd 2 copy/_ym xd/_xm xd tr n 1
lj 1 lc 1 w 0 8{dup}repeat sc fc bc}bd/clpg{n M dup neg 3 1 roll 0
l 0 exch l 0 l P clip n}bd/SB{_o 0 ne _d 0 ne and{gs []0 setdash _br
_bg _bb srgbc stroke gr}if}bd/S{_r _g _b srgbc SB stroke}bd/SX{_r _g
_b srgbc _cm currentmatrix pop _m setmatrix SB stroke _cm setmatrix}bd/d{/_d
xd setdash}bd/c{3{255 div 3 1 roll}repeat}bd/sc{c/_b xd/_g xd/_r xd}bd/bc{c/_bb
xd/_bg xd/_br xd}bd/fc{c/_B xd/_G xd/_R xd}bd/tc{c srgbc}bd/f{0 eq{eofill}{fill}ifelse}bd/D{{l}repeat}bd/PT{dup
statusdict exch known{cvx exec}{pop}ifelse}bd/TT{dup statusdict exch
known{cvx exec}{pop pop}ifelse}bd/SM{_m currentmatrix pop}bd/RM{_m
setmatrix}bd/sp{dup 0 eq{/_pb xd}{255 div/_gr xd/_pb xd 100 div/_f
xd}ifelse}bd/B{/_t2 xd/_t3 xd _pb _t2 _t3 8 idiv add get 1 7 _t3 8
mod sub bitshift and 0 ne}bd/SF{2{1 add 4 mul cvi exch}repeat B{0}{1}ifelse}bd/pb{_gr
setgray _f _a/SF load setscreen}bd
/s{100 div/_t xd[_t 0 0 _t 1 _t sub _xp 2 div mul 1 _t sub _yp 2 div
mul]concat}bd
/_am matrix def/a{_am currentmatrix pop/_t1 xd/_t2 xd/_t3 exch 10
div def/_t4 exch 10 div def tr _t2 _t1 neg scale/_t xd _t 1 eq{0 0
M}if 0 0 1 _t4 _t3 arc _t 0 ne{P}if _am setmatrix}bd/an{_am currentmatrix
pop/_t1 xd/_t2 xd/_t3 exch 10 div def/_t4 exch 10 div def tr _t2 _t1
neg scale/_t xd _t 1 eq{0 0 M}if 0 0 1 _t4 _t3 arcn _t 0 ne{P}if _am
setmatrix}bd
/_em matrix def/E{_em currentmatrix pop/_t1 xd/_t2 xd/_t3 xd/_t4 xd
_t4 _t3 _t1 sub M _t4 _t3 tr _t2 _t1 neg scale 0 0 1 90 450 arc _em
setmatrix}bd
/_rm matrix def/cn{_rm currentmatrix pop tr x3 y3 scale arc _rm setmatrix}bd/rr{2
div/y3 xd 2 div/x3 xd/y2 xd/x2 xd/y1 xd/x1 xd x1 y1 y3 add M 0 0 1
180 270 x1 x3 add y1 y3 add cn x2 x3 sub y1 L 0 0 1 270 0 x2 x3 sub
y1 y3 add cn x2 y2 y3 sub L 0 0 1 0 90 x2 x3 sub y2 y3 sub cn x1 x3
add y2 L 0 0 1 90 180 x1 x3 add y2 y3 sub cn}bd
/cp{/_t xd{{gs _t 0 eq{eoclip}{clip}ifelse}stopped{gr currentflat
1 add setflat}{exit}ifelse}loop n}bd/p{P dup 0 eq{S pop}{_R _G _B srgbc
_pb 0 eq{2 eq{gs f gr S}{f}ifelse}{_pb 1 eq{_o 1 eq{gs _br _bg _bb
srgbc 1 index f gr}if _R _G _B srgbc 2 eq{gs hf gr S}{hf}ifelse n}{2
eq{true pbc S}{false pbc n}ifelse}ifelse}ifelse}ifelse}bd/px{P dup
0 eq{SX pop}{_R _G _B srgbc _pb 0 eq{2 eq{gs f gr SX}{f}ifelse}{_pb
1 eq{_o 1 eq{gs _br _bg _bb srgbc 1 index f gr}if 2 eq{gs _m setmatrix
hf gr SX}{gs _m setmatrix hf gr}ifelse n}{2 eq{true pbc SX}{false pbc
n}ifelse}ifelse}ifelse}ifelse}bd
/spc{dup 0 eq{/_pb xd}{1 eq{255 div/_gr xd/_pb xd/_pg 0 def/_pr 0
def}{/_pb xd/_pg xd/_pr xd}ifelse 1 add/_pq xd}ifelse}bd/_im{[_t 0
0 _t 0 0]}bd/bbx{_pq 75 mul _rz div/_t xd pathbbox 1 add cvi 31 or/_ury
xd 1 add cvi 31 or/_urx xd 1 sub cvi -32 and/_lly xd 1 sub cvi -32
and/_llx xd/_dH _ury _lly sub def/_dW _urx _llx sub def}bd/xp{/_pbs
xd _row _pbs mul _pbs getinterval/_s xd 0 _pbs 8 _pbs div dup _sW add
1 sub exch div cvi dup/str exch string def _pbs sub{str exch _s putinterval}for
str}bd/pbc{/__t save def{_m setmatrix}if 0 eq{eoclip}{clip}ifelse/_row
0 def bbx _llx _lly tr _t _dW mul ceiling cvi 8 add dup 8 mod sub dup/_sW
xd _t _dH mul ceiling cvi 8 add dup 8 mod sub _pg 0 ne{8 _im{_pr 8
xp}{_pg 8 xp}{_pb 8 xp/_row _row 1 add 8 mod def}true 3 colorimage}{1
_o eq{gs _br _bg _bb srgbc 2 copy true _im{_pb 1 xp/_row _row 1 add
8 mod def}imagemask gr}if _R _G _B srgbc false _im{_pb 1 xp/_row _row
1 add 8 mod def}imagemask}ifelse __t restore}bd
/hp{dup/_hb xd -1 eq{/_pb 0 def}{1 add/_pq xd/_pb 1 def}ifelse}bd/vert{X0
_w X1{dup Y0 M Y1 L stroke}for}bd/horz{Y0 _w Y1{dup X0 exch M X1 exch
L stroke}for}bd/fdiag{X0 _w X1{Y0 M X1 X0 sub dup l stroke}for Y0 _w
Y1{X0 exch M Y1 Y0 sub dup l stroke}for}bd/bdiag{X0 _w X1{Y1 M X1 X0
sub dup neg l stroke}for Y0 _w Y1{X0 exch M Y1 Y0 sub dup neg l stroke}for}bd/AU{1
add cvi 31 or}bd/AD{1 sub cvi -32 and}bd/hf{pathbbox 1 add cvi 31 or/Y1
xd 1 add cvi 31 or/X1 xd 1 sub cvi -32 and/Y0 xd 1 sub cvi -32 and/X0
xd 2 w [] 0 setdash/_w _rz 20 div 8 div 2 mul _pq div round 8 mul def
cp _hb 0 eq{horz}if _hb 1 eq{vert}if _hb 2 eq{fdiag}if _hb 3 eq{bdiag}if
_hb 4 eq{horz vert}if _hb 5 eq{fdiag bdiag}if gr}bd
/sp{dup 0 eq{/_pb xd}{255 div/_gr xd/_pb xd 100 div/_f xd}ifelse}bd/B{/_t2
xd/_t3 xd _pb _t2 _t3 8 idiv add get 1 7 _t3 8 mod sub bitshift and
0 ne}bd/SF{2{1 add 4 mul cvi exch}repeat B{0}{1}ifelse}bd/pb{_gr setgray
_f _a/SF load setscreen}bd
/SCA{10 div dup dup/_a xd/_A xd/__na xd currentscreen/__p xd/__a xd/__f
xd __f __na/__p load setscreen}bd
/SCF{100 div dup dup/_F xd/__nf xd currentscreen/__p xd/__a xd/__f
xd __nf __a/__p load setscreen}bd
/makeOblique 0 def/bC 1 def/makeoutlinedict 7 dict def/mbf{makeoutlinedict
begin/uniqueid exch def/strokewidth exch def/newfontname exch def/basefontname
exch def/basefontdict basefontname findfont def/numentries basefontdict
maxlength 1 add def basefontdict/UniqueID known not{/numentries numentries
1 add def}if/outfontdict numentries dict def basefontdict{exch dup/FID
ne{exch outfontdict 3 1 roll put}{pop pop}ifelse}forall outfontdict/FontName
newfontname put outfontdict/PaintType 2 put outfontdict/StrokeWidth
strokewidth put outfontdict/UniqueID uniqueid put newfontname outfontdict
definefont pop end}bd/usf{findfont exch makeOblique 0 ne{/makeOblique
makeOblique sin neg def [1 0 makeOblique 1 0 0] matrix concatmatrix/makeOblique
0 def}if makefont setfont}bd/sf{/_nf xd FontDirectory _nf known{pop}{findfont
dup maxlength dict/_nfd xd{exch dup/FID ne{exch _nfd 3 1 roll put}{pop
pop}ifelse}forall _nfd dup/FontName _nf put/Encoding WE put _nf _nfd
definefont pop}ifelse _nf usf}bd/us{gs [] 0 setdash 0 setlinecap w
0 exch M _w _be add 0 l currentrgbcolor sC S gr}bd/ob{n/_h xd _cx _cy
M _w 0 l 0 _h l _w neg 0 l P tc fill}bd/st{dup/_be xd 0 ne{/_bc xd}if/_vt
xd dup stringwidth pop/_sw xd dup length dup 1 gt{1 sub}if/_sl xd _be
0 eq{_w _sw sub _sl div 0 _vt{exch}if 3 -1 roll ashow}{_be _bc div
0 _vt{exch}if 32 _w _sw sub _sl div 0 _vt{exch}if 6 -1 roll awidthshow}ifelse}bd/sb{gs
M 10 div rotate/_w xd currentpoint/_cy xd/_cx xd 0 ne{ob}if tc _cx
_cy 3 -1 roll add tr 0 0 M st dup 4 and 4 eq{4 sub/_dx _cx _px sub
def _dx neg 0 tr/_w _w _dx add def}if dup 0 eq{pop}{dup 1 eq{pop us}{2
eq{us}{dup 3 1 roll us us}ifelse}ifelse}ifelse _cx/_px xd gr}bd
/cl{gs n M dup neg 3 1 roll 0 l 0 exch l 0 l P clip n}bd/clx{gs n
_cm currentmatrix pop _m setmatrix n M dup neg 3 1 roll 0 l 0 exch
l 0 l P clip _cm setmatrix n}bd
/sbb{/_t save def _neg 0 eq{0}{1}ifelse setgray tr/_dW xd/_dH xd dup/_sW
xd 8 _bs div dup _sW add 1 sub exch div cvi/str exch string def exch
dup/_sH xd 3 -1 roll{/_im{[_sW _dW div 0 0 _sH neg _dH div 0 _sH]}bd}{/_im{[_sW
_dW div 0 0 _sH _dH div 0 0]}bd}ifelse _bs 1 eq{false}{_bs}ifelse _im{currentfile
str readhexstring pop}_bs 1 eq{imagemask}{image}ifelse _t restore}bd
/scbb{/_t save def tr/_dW xd/_dH xd dup/_sW xd 8 _bs div dup _sW add
1 sub exch div cvi 3 mul/str exch string def exch dup/_sH xd _bs 4
-1 roll{/_im{[_sW _dW div 0 0 _sH neg _dH div 0 _sH]}bd}{/_im{[_sW
_dW div 0 0 _sH _dH div 0 0]}bd}ifelse _im{currentfile str readhexstring
pop}false 3 colorimage _t restore}bd
/mr{2 copy -1 eq{_yp _ym 2 mul sub}{0}ifelse exch -1 eq{_xp _xm 2
mul sub}{0}ifelse exch tr scale}bd
/WE 256 array def StandardEncoding WE copy pop WE dup
0[/grave/acute/circumflex/tilde/macron/breve/dotaccent/dieresis/ring/cedilla
/hungarumlaut/ogonek/caron/dotlessi
]putinterval dup
39/quotesingle put dup  96/grave put dup 124/bar put dup
145[/quoteleft/quoteright/quotedblleft/quotedblright/bullet/endash/emdash]putinterval 
160[/space/exclamdown/cent/sterling/currency/yen/brokenbar/section/dieresis
/copyright/ordfeminine/guillemotleft/logicalnot/endash/registered/macron/ring
/plusminus/twosuperior/threesuperior/acute/mu/paragraph/periodcentered
/cedilla/onesuperior/ordmasculine/guillemotright/onequarter/onehalf
/threequarters/questiondown/Agrave/Aacute/Acircumflex/Atilde/Adieresis/Aring
/AE/Ccedilla/Egrave/Eacute/Ecircumflex/Edieresis/Igrave/Iacute/Icircumflex
/Idieresis/Eth/Ntilde/Ograve/Oacute/Ocircumflex/Otilde/Odieresis/circumflex
/Oslash/Ugrave/Uacute/Ucircumflex/Udieresis/Yacute/Thorn/germandbls/agrave
/aacute/acircumflex/atilde/adieresis/aring/ae/ccedilla/egrave/eacute
/ecircumflex/edieresis/igrave/iacute/icircumflex/idieresis/eth/ntilde/ograve
/oacute/ocircumflex/otilde/odieresis/tilde/oslash/ugrave/uacute/ucircumflex
/udieresis/yacute/thorn/ydieresis]putinterval
/TTEV 256 array def TTEV
0 [/G00/G01/G02/G03/G04/G05/G06/G07/G08/G09/G0a/G0b/G0c/G0d/G0e/G0f
   /G10/G11/G12/G13/G14/G15/G16/G17/G18/G19/G1a/G1b/G1c/G1d/G1e/G1f
   /G20/G21/G22/G23/G24/G25/G26/G27/G28/G29/G2a/G2b/G2c/G2d/G2e/G2f
   /G30/G31/G32/G33/G34/G35/G36/G37/G38/G39/G3a/G3b/G3c/G3d/G3e/G3f
   /G40/G41/G42/G43/G44/G45/G46/G47/G48/G49/G4a/G4b/G4c/G4d/G4e/G4f
   /G50/G51/G52/G53/G54/G55/G56/G57/G58/G59/G5a/G5b/G5c/G5d/G5e/G5f
   /G60/G61/G62/G63/G64/G65/G66/G67/G68/G69/G6a/G6b/G6c/G6d/G6e/G6f
   /G70/G71/G72/G73/G74/G75/G76/G77/G78/G79/G7a/G7b/G7c/G7d/G7e/G7f
   /G80/G81/G82/G83/G84/G85/G86/G87/G88/G89/G8a/G8b/G8c/G8d/G8e/G8f
   /G90/G91/G92/G93/G94/G95/G96/G97/G98/G99/G9a/G9b/G9c/G9d/G9e/G9f
   /Ga0/Ga1/Ga2/Ga3/Ga4/Ga5/Ga6/Ga7/Ga8/Ga9/Gaa/Gab/Gac/Gad/Gae/Gaf
   /Gb0/Gb1/Gb2/Gb3/Gb4/Gb5/Gb6/Gb7/Gb8/Gb9/Gba/Gbb/Gbc/Gbd/Gbe/Gbf
   /Gc0/Gc1/Gc2/Gc3/Gc4/Gc5/Gc6/Gc7/Gc8/Gc9/Gca/Gcb/Gcc/Gcd/Gce/Gcf
   /Gd0/Gd1/Gd2/Gd3/Gd4/Gd5/Gd6/Gd7/Gd8/Gd9/Gda/Gdb/Gdc/Gdd/Gde/Gdf
   /Ge0/Ge1/Ge2/Ge3/Ge4/Ge5/Ge6/Ge7/Ge8/Ge9/Gea/Geb/Gec/Ged/Gee/Gef
   /Gf0/Gf1/Gf2/Gf3/Gf4/Gf5/Gf6/Gf7/Gf8/Gf9/Gfa/Gfb/Gfc/Gfd/Gfe/Gff
] putinterval
end
MGXPS_2.8_EPS begin
0 0 1 0 3307 4677 400 -1 0 0 1 Cne
4551 3134 0 0 clpg
0 0 0 sc
/_o 0 def
255 255 255 bc
[5 13 ]0 2 d
156 2 M 218 23 272 57 299 85 bz
S
157 480 M 115 486 57 493 0 483 bz
S
[]0 0 d
300 86 M 378 174 374 208 368 302 bz
359 398 269 455 158 480 bz
S
402 0 M 0 0 -55 495 2 D
S
197 283 3 3 E
0 0 0 fc
0 2 p
198 284 M 0 0 170 18 2 D
S
198 283 M 0 0 146 -139 2 D
S
[5 13 ]0 2 d
344 144 M 0 0 47 -46 2 D
S
[]0 0 d
367 301 3 3 E
0 2 p
330 144 M 0 0 14 0 2 D
S
344 144 M 0 0 0 13 2 D
S
208 295 M 0 0 -8 -11 2 D
S
200 284 M 0 0 10 -8 2 D
S
215 281 M 0 0 -14 0 2 D
S
201 281 M 0 0 0 -13 2 D
S
357 291 M 0 0 8 11 2 D
S
365 302 M 0 0 -10 8 2 D
S
224 260 M 230 267 235 274 232 287 bz
S
end
EndOfTheIncludedPostscriptMagicCookie

\closepsdump

\psdump{hyper.eps}
statusdict begin/waittimeout 0 def end{/MGXPS_2.8_EPS begin}stopped{/MGXPS_2.8_EPS
200 dict def}if MGXPS_2.8_EPS begin/bd{bind def}bind def/xd{exch def}bd/ld{load
def}bd/M/moveto ld/m/rmoveto ld/L/lineto ld/l/rlineto ld/w/setlinewidth
ld/n/newpath ld/P/closepath ld/tr/translate ld/gs/gsave ld/gr/grestore
ld/lc/setlinecap ld/lj/setlinejoin ld/ml/setmiterlimit ld/TM/concat
ld/F/fill ld/bz/curveto ld/_m matrix def/_cm matrix def/_d 0 def/_o
0 def/_pb 0 def/_hb -1 def currentscreen/_P xd dup/_A xd/__A xd dup/_F
xd/__F xd/srgbc{_neg 0 ne{3{neg 1 add 3 1 roll}repeat}if setrgbcolor}bd/ssc{gs
0 0 M 0 _yp neg _xp 0 0 _yp 3 D closepath fill gr 0 0 0 srgbc}bd/Cne{/_bs
xd/_neg where{/_neg get/_neg xd pop}{pop 0/_neg xd}ifelse dup 0 ne{setflat}{pop}ifelse
dup -1 eq{pop}{pop pop}ifelse/_rz xd/_yp xd/_xp xd/_a 0 def 0 ne{/_a
90 def -90 rotate _xp _rz div 72 mul neg 0 tr}if 0 _yp _rz div 72 mul
tr 72 _rz div dup neg scale/#copies xd 2 copy/_ym xd/_xm xd tr n 1
lj 1 lc 1 w 0 8{dup}repeat sc fc bc}bd/clpg{n M dup neg 3 1 roll 0
l 0 exch l 0 l P clip n}bd/SB{_o 0 ne _d 0 ne and{gs []0 setdash _br
_bg _bb srgbc stroke gr}if}bd/S{_r _g _b srgbc SB stroke}bd/SX{_r _g
_b srgbc _cm currentmatrix pop _m setmatrix SB stroke _cm setmatrix}bd/d{/_d
xd setdash}bd/c{3{255 div 3 1 roll}repeat}bd/sc{c/_b xd/_g xd/_r xd}bd/bc{c/_bb
xd/_bg xd/_br xd}bd/fc{c/_B xd/_G xd/_R xd}bd/tc{c srgbc}bd/f{0 eq{eofill}{fill}ifelse}bd/D{{l}repeat}bd/PT{dup
statusdict exch known{cvx exec}{pop}ifelse}bd/TT{dup statusdict exch
known{cvx exec}{pop pop}ifelse}bd/SM{_m currentmatrix pop}bd/RM{_m
setmatrix}bd/sp{dup 0 eq{/_pb xd}{255 div/_gr xd/_pb xd 100 div/_f
xd}ifelse}bd/B{/_t2 xd/_t3 xd _pb _t2 _t3 8 idiv add get 1 7 _t3 8
mod sub bitshift and 0 ne}bd/SF{2{1 add 4 mul cvi exch}repeat B{0}{1}ifelse}bd/pb{_gr
setgray _f _a/SF load setscreen}bd
/s{100 div/_t xd[_t 0 0 _t 1 _t sub _xp 2 div mul 1 _t sub _yp 2 div
mul]concat}bd
/_am matrix def/a{_am currentmatrix pop/_t1 xd/_t2 xd/_t3 exch 10
div def/_t4 exch 10 div def tr _t2 _t1 neg scale/_t xd _t 1 eq{0 0
M}if 0 0 1 _t4 _t3 arc _t 0 ne{P}if _am setmatrix}bd/an{_am currentmatrix
pop/_t1 xd/_t2 xd/_t3 exch 10 div def/_t4 exch 10 div def tr _t2 _t1
neg scale/_t xd _t 1 eq{0 0 M}if 0 0 1 _t4 _t3 arcn _t 0 ne{P}if _am
setmatrix}bd
/_em matrix def/E{_em currentmatrix pop/_t1 xd/_t2 xd/_t3 xd/_t4 xd
_t4 _t3 _t1 sub M _t4 _t3 tr _t2 _t1 neg scale 0 0 1 90 450 arc _em
setmatrix}bd
/_rm matrix def/cn{_rm currentmatrix pop tr x3 y3 scale arc _rm setmatrix}bd/rr{2
div/y3 xd 2 div/x3 xd/y2 xd/x2 xd/y1 xd/x1 xd x1 y1 y3 add M 0 0 1
180 270 x1 x3 add y1 y3 add cn x2 x3 sub y1 L 0 0 1 270 0 x2 x3 sub
y1 y3 add cn x2 y2 y3 sub L 0 0 1 0 90 x2 x3 sub y2 y3 sub cn x1 x3
add y2 L 0 0 1 90 180 x1 x3 add y2 y3 sub cn}bd
/cp{/_t xd{{gs _t 0 eq{eoclip}{clip}ifelse}stopped{gr currentflat
1 add setflat}{exit}ifelse}loop n}bd/p{P dup 0 eq{S pop}{_R _G _B srgbc
_pb 0 eq{2 eq{gs f gr S}{f}ifelse}{_pb 1 eq{_o 1 eq{gs _br _bg _bb
srgbc 1 index f gr}if _R _G _B srgbc 2 eq{gs hf gr S}{hf}ifelse n}{2
eq{true pbc S}{false pbc n}ifelse}ifelse}ifelse}ifelse}bd/px{P dup
0 eq{SX pop}{_R _G _B srgbc _pb 0 eq{2 eq{gs f gr SX}{f}ifelse}{_pb
1 eq{_o 1 eq{gs _br _bg _bb srgbc 1 index f gr}if 2 eq{gs _m setmatrix
hf gr SX}{gs _m setmatrix hf gr}ifelse n}{2 eq{true pbc SX}{false pbc
n}ifelse}ifelse}ifelse}ifelse}bd
/spc{dup 0 eq{/_pb xd}{1 eq{255 div/_gr xd/_pb xd/_pg 0 def/_pr 0
def}{/_pb xd/_pg xd/_pr xd}ifelse 1 add/_pq xd}ifelse}bd/_im{[_t 0
0 _t 0 0]}bd/bbx{_pq 75 mul _rz div/_t xd pathbbox 1 add cvi 31 or/_ury
xd 1 add cvi 31 or/_urx xd 1 sub cvi -32 and/_lly xd 1 sub cvi -32
and/_llx xd/_dH _ury _lly sub def/_dW _urx _llx sub def}bd/xp{/_pbs
xd _row _pbs mul _pbs getinterval/_s xd 0 _pbs 8 _pbs div dup _sW add
1 sub exch div cvi dup/str exch string def _pbs sub{str exch _s putinterval}for
str}bd/pbc{/__t save def{_m setmatrix}if 0 eq{eoclip}{clip}ifelse/_row
0 def bbx _llx _lly tr _t _dW mul ceiling cvi 8 add dup 8 mod sub dup/_sW
xd _t _dH mul ceiling cvi 8 add dup 8 mod sub _pg 0 ne{8 _im{_pr 8
xp}{_pg 8 xp}{_pb 8 xp/_row _row 1 add 8 mod def}true 3 colorimage}{1
_o eq{gs _br _bg _bb srgbc 2 copy true _im{_pb 1 xp/_row _row 1 add
8 mod def}imagemask gr}if _R _G _B srgbc false _im{_pb 1 xp/_row _row
1 add 8 mod def}imagemask}ifelse __t restore}bd
/hp{dup/_hb xd -1 eq{/_pb 0 def}{1 add/_pq xd/_pb 1 def}ifelse}bd/vert{X0
_w X1{dup Y0 M Y1 L stroke}for}bd/horz{Y0 _w Y1{dup X0 exch M X1 exch
L stroke}for}bd/fdiag{X0 _w X1{Y0 M X1 X0 sub dup l stroke}for Y0 _w
Y1{X0 exch M Y1 Y0 sub dup l stroke}for}bd/bdiag{X0 _w X1{Y1 M X1 X0
sub dup neg l stroke}for Y0 _w Y1{X0 exch M Y1 Y0 sub dup neg l stroke}for}bd/AU{1
add cvi 31 or}bd/AD{1 sub cvi -32 and}bd/hf{pathbbox 1 add cvi 31 or/Y1
xd 1 add cvi 31 or/X1 xd 1 sub cvi -32 and/Y0 xd 1 sub cvi -32 and/X0
xd 2 w [] 0 setdash/_w _rz 20 div 8 div 2 mul _pq div round 8 mul def
cp _hb 0 eq{horz}if _hb 1 eq{vert}if _hb 2 eq{fdiag}if _hb 3 eq{bdiag}if
_hb 4 eq{horz vert}if _hb 5 eq{fdiag bdiag}if gr}bd
/sp{dup 0 eq{/_pb xd}{255 div/_gr xd/_pb xd 100 div/_f xd}ifelse}bd/B{/_t2
xd/_t3 xd _pb _t2 _t3 8 idiv add get 1 7 _t3 8 mod sub bitshift and
0 ne}bd/SF{2{1 add 4 mul cvi exch}repeat B{0}{1}ifelse}bd/pb{_gr setgray
_f _a/SF load setscreen}bd
/SCA{10 div dup dup/_a xd/_A xd/__na xd currentscreen/__p xd/__a xd/__f
xd __f __na/__p load setscreen}bd
/SCF{100 div dup dup/_F xd/__nf xd currentscreen/__p xd/__a xd/__f
xd __nf __a/__p load setscreen}bd
/makeOblique 0 def/bC 1 def/makeoutlinedict 7 dict def/mbf{makeoutlinedict
begin/uniqueid exch def/strokewidth exch def/newfontname exch def/basefontname
exch def/basefontdict basefontname findfont def/numentries basefontdict
maxlength 1 add def basefontdict/UniqueID known not{/numentries numentries
1 add def}if/outfontdict numentries dict def basefontdict{exch dup/FID
ne{exch outfontdict 3 1 roll put}{pop pop}ifelse}forall outfontdict/FontName
newfontname put outfontdict/PaintType 2 put outfontdict/StrokeWidth
strokewidth put outfontdict/UniqueID uniqueid put newfontname outfontdict
definefont pop end}bd/usf{findfont exch makeOblique 0 ne{/makeOblique
makeOblique sin neg def [1 0 makeOblique 1 0 0] matrix concatmatrix/makeOblique
0 def}if makefont setfont}bd/sf{/_nf xd FontDirectory _nf known{pop}{findfont
dup maxlength dict/_nfd xd{exch dup/FID ne{exch _nfd 3 1 roll put}{pop
pop}ifelse}forall _nfd dup/FontName _nf put/Encoding WE put _nf _nfd
definefont pop}ifelse _nf usf}bd/us{gs [] 0 setdash 0 setlinecap w
0 exch M _w _be add 0 l currentrgbcolor sC S gr}bd/ob{n/_h xd _cx _cy
M _w 0 l 0 _h l _w neg 0 l P tc fill}bd/st{dup/_be xd 0 ne{/_bc xd}if/_vt
xd dup stringwidth pop/_sw xd dup length dup 1 gt{1 sub}if/_sl xd _be
0 eq{_w _sw sub _sl div 0 _vt{exch}if 3 -1 roll ashow}{_be _bc div
0 _vt{exch}if 32 _w _sw sub _sl div 0 _vt{exch}if 6 -1 roll awidthshow}ifelse}bd/sb{gs
M 10 div rotate/_w xd currentpoint/_cy xd/_cx xd 0 ne{ob}if tc _cx
_cy 3 -1 roll add tr 0 0 M st dup 4 and 4 eq{4 sub/_dx _cx _px sub
def _dx neg 0 tr/_w _w _dx add def}if dup 0 eq{pop}{dup 1 eq{pop us}{2
eq{us}{dup 3 1 roll us us}ifelse}ifelse}ifelse _cx/_px xd gr}bd
/cl{gs n M dup neg 3 1 roll 0 l 0 exch l 0 l P clip n}bd/clx{gs n
_cm currentmatrix pop _m setmatrix n M dup neg 3 1 roll 0 l 0 exch
l 0 l P clip _cm setmatrix n}bd
/sbb{/_t save def _neg 0 eq{0}{1}ifelse setgray tr/_dW xd/_dH xd dup/_sW
xd 8 _bs div dup _sW add 1 sub exch div cvi/str exch string def exch
dup/_sH xd 3 -1 roll{/_im{[_sW _dW div 0 0 _sH neg _dH div 0 _sH]}bd}{/_im{[_sW
_dW div 0 0 _sH _dH div 0 0]}bd}ifelse _bs 1 eq{false}{_bs}ifelse _im{currentfile
str readhexstring pop}_bs 1 eq{imagemask}{image}ifelse _t restore}bd
/scbb{/_t save def tr/_dW xd/_dH xd dup/_sW xd 8 _bs div dup _sW add
1 sub exch div cvi 3 mul/str exch string def exch dup/_sH xd _bs 4
-1 roll{/_im{[_sW _dW div 0 0 _sH neg _dH div 0 _sH]}bd}{/_im{[_sW
_dW div 0 0 _sH _dH div 0 0]}bd}ifelse _im{currentfile str readhexstring
pop}false 3 colorimage _t restore}bd
/mr{2 copy -1 eq{_yp _ym 2 mul sub}{0}ifelse exch -1 eq{_xp _xm 2
mul sub}{0}ifelse exch tr scale}bd
/WE 256 array def StandardEncoding WE copy pop WE dup
0[/grave/acute/circumflex/tilde/macron/breve/dotaccent/dieresis/ring/cedilla
/hungarumlaut/ogonek/caron/dotlessi
]putinterval dup
39/quotesingle put dup  96/grave put dup 124/bar put dup
145[/quoteleft/quoteright/quotedblleft/quotedblright/bullet/endash/emdash]putinterval 
160[/space/exclamdown/cent/sterling/currency/yen/brokenbar/section/dieresis
/copyright/ordfeminine/guillemotleft/logicalnot/endash/registered/macron/ring
/plusminus/twosuperior/threesuperior/acute/mu/paragraph/periodcentered
/cedilla/onesuperior/ordmasculine/guillemotright/onequarter/onehalf
/threequarters/questiondown/Agrave/Aacute/Acircumflex/Atilde/Adieresis/Aring
/AE/Ccedilla/Egrave/Eacute/Ecircumflex/Edieresis/Igrave/Iacute/Icircumflex
/Idieresis/Eth/Ntilde/Ograve/Oacute/Ocircumflex/Otilde/Odieresis/circumflex
/Oslash/Ugrave/Uacute/Ucircumflex/Udieresis/Yacute/Thorn/germandbls/agrave
/aacute/acircumflex/atilde/adieresis/aring/ae/ccedilla/egrave/eacute
/ecircumflex/edieresis/igrave/iacute/icircumflex/idieresis/eth/ntilde/ograve
/oacute/ocircumflex/otilde/odieresis/tilde/oslash/ugrave/uacute/ucircumflex
/udieresis/yacute/thorn/ydieresis]putinterval
/TTEV 256 array def TTEV
0 [/G00/G01/G02/G03/G04/G05/G06/G07/G08/G09/G0a/G0b/G0c/G0d/G0e/G0f
   /G10/G11/G12/G13/G14/G15/G16/G17/G18/G19/G1a/G1b/G1c/G1d/G1e/G1f
   /G20/G21/G22/G23/G24/G25/G26/G27/G28/G29/G2a/G2b/G2c/G2d/G2e/G2f
   /G30/G31/G32/G33/G34/G35/G36/G37/G38/G39/G3a/G3b/G3c/G3d/G3e/G3f
   /G40/G41/G42/G43/G44/G45/G46/G47/G48/G49/G4a/G4b/G4c/G4d/G4e/G4f
   /G50/G51/G52/G53/G54/G55/G56/G57/G58/G59/G5a/G5b/G5c/G5d/G5e/G5f
   /G60/G61/G62/G63/G64/G65/G66/G67/G68/G69/G6a/G6b/G6c/G6d/G6e/G6f
   /G70/G71/G72/G73/G74/G75/G76/G77/G78/G79/G7a/G7b/G7c/G7d/G7e/G7f
   /G80/G81/G82/G83/G84/G85/G86/G87/G88/G89/G8a/G8b/G8c/G8d/G8e/G8f
   /G90/G91/G92/G93/G94/G95/G96/G97/G98/G99/G9a/G9b/G9c/G9d/G9e/G9f
   /Ga0/Ga1/Ga2/Ga3/Ga4/Ga5/Ga6/Ga7/Ga8/Ga9/Gaa/Gab/Gac/Gad/Gae/Gaf
   /Gb0/Gb1/Gb2/Gb3/Gb4/Gb5/Gb6/Gb7/Gb8/Gb9/Gba/Gbb/Gbc/Gbd/Gbe/Gbf
   /Gc0/Gc1/Gc2/Gc3/Gc4/Gc5/Gc6/Gc7/Gc8/Gc9/Gca/Gcb/Gcc/Gcd/Gce/Gcf
   /Gd0/Gd1/Gd2/Gd3/Gd4/Gd5/Gd6/Gd7/Gd8/Gd9/Gda/Gdb/Gdc/Gdd/Gde/Gdf
   /Ge0/Ge1/Ge2/Ge3/Ge4/Ge5/Ge6/Ge7/Ge8/Ge9/Gea/Geb/Gec/Ged/Gee/Gef
   /Gf0/Gf1/Gf2/Gf3/Gf4/Gf5/Gf6/Gf7/Gf8/Gf9/Gfa/Gfb/Gfc/Gfd/Gfe/Gff
] putinterval
end
MGXPS_2.8_EPS begin
0 0 1 0 3307 4677 400 -1 0 0 1 Cne
4551 3134 0 0 clpg
0 0 0 sc
/_o 0 def
255 255 255 bc
[5 13 ]0 2 d
156 2 M 218 23 272 57 299 85 bz
S
157 480 M 115 486 57 493 0 483 bz
S
[]0 0 d
300 86 M 378 174 374 208 368 302 bz
359 398 269 455 158 480 bz
S
402 0 M 0 0 -55 495 2 D
S
197 283 3 3 E
0 0 0 fc
0 2 p
198 284 M 0 0 170 18 2 D
S
367 301 3 3 E
0 2 p
199 280 169 169 E
0 p
132 51 3 3 E
0 2 p
133 52 M 0 0 236 250 2 D
S
198 284 M 0 0 82 -77 2 D
S
278 208 3 3 E
0 2 p
end
EndOfTheIncludedPostscriptMagicCookie

\closepsdump

\usepackage[dvips]{epsfig}
\usepackage{amstex,ifthen,amssymb,calc,epic}

\allowdisplaybreaks

\makeatletter
\renewcommand{\@seccntformat}[1]{\csname the#1\endcsname.\hspace{1em}}
\renewcommand\section{\@startsection {section}{1}{\z@}
                                   {-3.5ex \@plus -1ex \@minus -.2ex}
                                   {2.3ex \@plus.2ex}
                                   {\normalfont\normalsize\scshape\centering}}
\renewcommand\subsection{\@startsection{subsection}{2}{\z@}
                                     {-3.25ex\@plus -1ex \@minus -.2ex}
                                     {1.5ex \@plus .2ex}
                                     {\normalfont\normalsize\scshape\centering}}
\renewcommand{\l@section}[2]{
  \ifnum \c@tocdepth >\z@
    \addpenalty{\@secpenalty}
    \addvspace{\p@}
    \setlength\@tempdima{1.5em}
    \begingroup
      \parindent \z@ \rightskip \@pnumwidth
      \parfillskip -\@pnumwidth
      \leavevmode \scshape
      \advance\leftskip\@tempdima
      \hskip -\leftskip
      #1\nobreak\hfil \nobreak\hbox to\@pnumwidth{\hss #2}\par
    \endgroup
  \fi}
\def\@maketitle{
  \newpage
  \null
  \vskip 2em
  \begin{center}
    {\sffamily \large \@title \par}
    \vskip 1.5em
      {\normalsize
      \lineskip .5em
        {\scshape \normalsize \@author}}
  \end{center}
  \par
  \vskip 1.5em}
\makeatother

\newcommand\poly{(-\Delta)^m\dir}%
\newcommand\real{\mathbb R}%
\newcommand\dup{\mathrm{d}}%
\newcommand\dir{|^{}_{\operatorname{DIR}}}%
\newcommand\power[1][\Omega]{H_{#1,m}}%
\newcommand\tc{\operatorname{tr}}%
\newcommand\ip[2]{\langle #1,#2\rangle}%
\newcommand\slref[1]{\textsl{\ref{#1}}}%
\newcommand\slcite[1]{\textsl{\cite{#1}}}%
\newcommand\finbox{~\hfill$\Box$}%
\newcommand\lowearlybox{\tag*{$\frac[0pt]{\displaystyle{\phantom{\big|}}}{%
  \displaystyle{\Box}}$}}%
\newcommand\bip[2]{\big\langle #1,#2\big\rangle}%
\newcommand\inradius{\operatorname{Inradius}}%
\newcommand\earlybox{\tag*{$\Box$}}%
\newcommand\spec{\operatorname{Spec}}%
\newcommand\integer{\mathbb Z}%
\renewcommand\eqref[1]{(\ref{#1})}%
\renewcommand\natural{\mathbb N}%
\newcommand\lin{\operatorname{lin}}%

\newcommand\afigwidth0%
\newcommand\afigheight0%
\newcommand\bfigwidth0%
\newcommand\bfigheight0%
\newlength\aeqspace%
\newlength\beqspace%
\newlength\border%
\newlength\tempeqspace%
\newcounter{tempactr}%
\newcounter{tempbctr}%
\newcounter{tempcctr}%
1%

\newcommand{\fig}[6]{%
  \setcounter{tempactr}{#6}%
  \setcounter{tempbctr}{200+\border/\unitlength}%
  \setcounter{tempcctr}{#6+\border/\unitlength}%
  \setlength\tempeqspace{#5}%
  \vskip\tempeqspace%
  \vskip\border%
  \centering%
  \begin{picture}(200,\value{tempactr})%
    \put(0,0){\epsfig{figure=#1.eps,%
                      width=200\unitlength}}%
    #3%
  \end{picture}%
  \vskip\border%
  \vskip\tempeqspace%
  \vskip\abovecaptionskip%
  \refstepcounter{figure}#2%
  \makebox[0pt][c]{\parbox[t]{1.5\textwidth}{\centering%
    \ifthenelse{\equal{#4}{}}{Figure \thefigure}{Figure \thefigure: #4}}}%
  \vskip\belowcaptionskip}%

\newcommand{\tab}[6]{%
  \setcounter{tempactr}{#6}%
  \setcounter{tempbctr}{200+\border/\unitlength}%
  \setcounter{tempcctr}{#6+\border/\unitlength}%
  \setlength\tempeqspace{#5}%
  \vskip\tempeqspace%
  \vskip\border%
  \centering%
  {\small%
  \begin{picture}(200,\value{tempactr})%
    \put(0,0){\begin{tabular}[b]{#2}%
      #3%
    \end{tabular}}%
  \end{picture}}%
  \vskip\border%
  \vskip\tempeqspace%
  \vskip\abovecaptionskip%
  \refstepcounter{table}#1%
  \makebox[0pt][c]{\parbox[t]{1.5\textwidth}{\centering%
    \ifthenelse{\equal{#4}{}}{Table \thetable}{Table \thetable: #4}}}%
  \vskip\belowcaptionskip}%

\newcommand{\dbledgram}[3]{{\setlength\border{0pt}
  \renewcommand\baselinestretch1\normalsize%
  \setlength\aeqspace{0pt}%
  \setlength\beqspace{0pt}%
  #1%
  \ifthenelse{\bfigheight > \afigheight}{%
    \setlength\aeqspace{\textwidth/200*(\bfigheight-\afigheight)}}{%
    \setlength\beqspace{\textwidth/200*(\afigheight-\bfigheight)}}%
  \begin{figure}[t]%
    \par%
    \noindent%
    \hfill%
    \setlength\unitlength{\textwidth/100*\afigwidth/200}%
    \begin{minipage}[t]{\textwidth/100*\afigwidth+2\border}%
      #2\aeqspace{200*\afigheight/\afigwidth}%
    \end{minipage}%
    \hfill%
    \setlength\unitlength{\textwidth/100*\bfigwidth/200}%
    \begin{minipage}[t]{\textwidth/100*\bfigwidth+2\border}%
      #3\beqspace{200*\bfigheight/\bfigwidth}%
    \end{minipage}%
    \hfill%
    \null%
  \end{figure}}}%

\newcounter{theorem}%
\newenvironment{thm}[1]{%
    \refstepcounter{theorem}%
    \textbf{#1 \thetheorem:}{}%
    \begin{itshape}}%
  {\end{itshape}}%

\newenvironment{proof}{\textbf{Proof:}{}}{}%

\newenvironment{note}{%
  \refstepcounter{theorem}%
  \textbf{Note \thetheorem:}{}}%
  {}%

\newenvironment{definition}[1][]{%
  \refstepcounter{theorem}%
  \ifthenelse{\equal{#1}{}}{\textbf{Definition \thetheorem:}{} }%
    {\textbf{Definition \thetheorem: (#1)}{} }%
  \begin{itshape}}%
  {\end{itshape}}%

\begin{document}

\title{The Hardy-Rellich Inequality for Polyharmonic Operators}
\author{M. P. Owen\\[1em]
  \slshape{\small Department of Mathematics, King's College London, Strand,
    London, WC2R 2LS}}
\maketitle

\begin{abstract}
  The Hardy-Rellich inequality given here generalizes a Hardy inequality of 
  Davies~\cite{Davi:1984}, from the case of the Dirichlet Laplacian of a 
  region $\Omega\subseteq\real^N$ to that of the higher order polyharmonic 
  operators with Dirichlet boundary conditions. The inequality yields some 
  immediate spectral information for the polyharmonic operators and also
  bounds on the trace of the associated semigroups and resolvents.
\end{abstract}%
\section{Introduction}%
The Hardy inequality originated in 1920 in~\cite{Hard:1920} as an integral 
inequality for functions defined on the real half-line. Its original 
representation can be easily reformulated, for $1<p<\infty$, as
\begin{equation}\label{eqn:orighardy}
  \int_0^\infty\frac{|f(x)|^p}{x^p}\dup x\le\left(\frac{p}{p-1}\right)^p\int_0^
  \infty|f'(x)|^p\dup x
\end{equation}
for all $f\in C_c^\infty((0,\infty))$. Since its appearance, various 
generalizations of particular aspects of the inequality have been made. 
In~\cite{Opic+Kufn:1990}, for example, there is a detailed treatment of 
weighted Hardy-type inequalities in an $L^p$ setting.

The Rellich inequality 
appeared first in~\cite{Rell:1956} as a generalization of 
inequality~\eqref{eqn:orighardy} to two derivatives. The simplest form of such 
an inequality is
\begin{equation}\label{eqn:origrellich}
  \int_0^\infty\frac{|f(x)|^p}{x^{2p}}\dup x\le\left(\frac{p^2}{(p-1)(2p-1)}
  \right)^p\int_0^\infty|f''(x)|^p\dup x
\end{equation}
for all $f\in C_c^\infty((0,\infty))$.

In this paper we study a generalization for all derivatives, within 
the $L^2$ setting. The variable $x$ in the denominator of 
inequalities~\eqref{eqn:orighardy} and~\eqref{eqn:origrellich} is replaced by 
a pseudodistance $a_m(x)$, where $m$ is the number of derivatives 
in the dominating integrand. We formulate our result as a Hardy-Rellich 
operator 
inequality, and use it as a tool in the spectral analysis of 
polyharmonic operators. The Rellich inequalities found in~\cite{Davi+Hinz:A1} 
concern a distinct but related class of operator inequalities.

In order to state our result properly, we need the following definitions, in 
which $\Omega$ denotes an open subset of $\real^N$:

\begin{definition}
  Let $Q_m$ be the closure of the quadratic form defined on $C_c^\infty(\Omega)
  \subseteq L^2(\Omega)$ by
  \[ Q_m(f)=\ip{(-\Delta)^mf}f.  \]
  The domain of the closure is the Sobolev space $W_0^{m,2}(\Omega)$. The 
  polyharmonic operator $\poly$ is defined as the non-negative self-adjoint 
  operator associated with $Q_m$. See~\slcite{Davi:1995} for details. Where
  the implied region is not contextually evident, the operator is denoted by 
  $\power$.
\end{definition}

The boundary conditions classically associated with the operators 
$\poly$ and $(-\Delta\dir)^m$ are different. The inequality
\begin{equation}\label{eqn:polylarger}
  \power^{}\ge H_{\Omega,1}^m
\end{equation}
may be verified by considering quadratic form domains.

\begin{definition}\label{def:pseud}
  Let $\omega\in S^{N-1}$ and define
  $d_\omega:\real^N\rightarrow(0,+\infty]$ by
  \begin{equation}\label{eqn:defd}
    d_\omega(x):=\min\{|s|:x+s\omega\not\in\Omega\}.
  \end{equation}
  Define the pseudodistances $a_m:\real^N\rightarrow(0,+\infty]$ for $1\le m
  \in\real$ by 
  \begin{equation}\label{eqn:pseud}
    a_m(x)=\left[\int_{S^{N-1}}d_\omega(x)^{-2m}\dup^{N-1}\omega\right]^{-1/2m
    },
  \end{equation}
  where $\dup^{N-1}\omega$ is the normalized surface measure on the unit 
  spherical shell $S^{N-1}$. 
\end{definition}

Davies~\cite{Davi:1984} proves the operator inequality 
\begin{equation}\label{eqn:oldhardy}
  (-\Delta)\dir\ge\frac N{4a_1^2}
\end{equation}
in the quadratic form sense, thus comparing the Dirichlet Laplacian with a 
multiplication operator which is large near the boundary of the region. In 
Theorem~\ref{thm:hardy} we generalize this to a Hardy-Rellich inequality
\begin{equation}\label{eqn:hardrell}
  \poly\ge\frac{(N+2m-2)(N+2m-4)\dots N(2m-1)(2m-3)\dots1}{4^ma_m^{2m}}
\end{equation}
for the polyharmonic operators $\poly=\power$ acting in $L^2(\Omega)$. As a 
special case, if $\Omega$ is convex we see that
\begin{equation}\label{eqn:intrconv}
  \poly\ge\frac{(2m-1)^2(2m-3)^2\dots1^2}{4^md^{2m}},
\end{equation}
where
\[ d(x):=\min\{|y-x|:y\not\in\Omega\}. \]
More generally, for regular regions, where the pseudodistances $a_m$ are 
comparable with the distance $d$, a similar inequality is valid. The constants 
in inequalities~\eqref{eqn:hardrell} and~\eqref{eqn:intrconv} are shown to be 
optimal.

In both~\cite{Davi:1985} and~\cite[section 1.9]{Davi:1989}, Davies uses the 
Hardy inequality~\eqref{eqn:oldhardy} to find an upper bound on the trace of 
the semigroup $e^{-H_{\Omega,1}t}$. Using the technique of decomposing a 
region with finite inradius into dyadic cubes he finds a similar lower bound. 
More explicitly,
\begin{equation}\label{eqn:oldtrace}
  (8\pi t)^{-N/2}\int_\Omega \exp[{-8\pi^2N^2t/d^2}]\le\tc[e^{-H_{\Omega,1}t}]
  \le(2\pi t)^{-N/2}\int_\Omega \exp[{-Nt/8a_1^2}].
\end{equation}
For regular regions $\Omega$ this yields an immediate equivalent condition for 
$\tc[e^{-H_{\Omega,1}t}]$ to be finite, and, as a corollary, a condition for 
$\tc[H_{\Omega,1}^{-\gamma}]$ to be finite.

In Sections~\ref{sec:lwbd} and~\ref{sec:upbd} we generalize 
inequality~\eqref{eqn:oldtrace} with 
some restrictions (see Condition~\ref{thm:bdker} and dependent results), to 
\begin{equation}\label{eqn:mytrace}
  b_{m,N}t^{-N/2m}\int_\Omega \exp[{-c_{m,N}td^{-2m}}]\le\tc[e^{-\power
  t}]\le b_{m,N}'t^{-N/{2m}}\int_\Omega \exp[{-c_{N,m}'ta_m^{-2m}}],
\end{equation}
where $b_{m,N}$, $c_{m,N}$, $b_{m,N}'$ and $c_{m,N}'$ are positive constants. 
This yields an equivalent condition for finite trace of both $e^{-\power t}$ 
and $\power^{-\gamma}$.%
\subsection*{Acknowledgments}
I wish to thank E B Davies for suggesting this problem and for his invaluable 
guidance and support during my research. I thank also Mark Ashbaugh for a 
useful comment. This research was funded by an EPSRC studentship.
\section{The Hardy-Rellich Operator Inequality}%
Our starting point is a one dimensional version of the Hardy-Rellich 
inequality in the $L^2$ setting. For $m=1$ and $m=2$, the following lemma 
respectively resembles inequalities~\eqref{eqn:orighardy} 
and~\eqref{eqn:origrellich}, where we set $p=2$.

\begin{thm}{Lemma}\label{thm:onedimhardy}
  Let $\Omega$ be an open (not necessarily connected) set in $\real$. Then
  \begin{equation}\label{eqn:onedimhardy}
    \frac{(2m-1)^2(2m-3)^2\dots1^2}{4^m}\int_\Omega\frac{|f(x)|^2}{d(x)^{2m}}
    \dup x\le\int_\Omega|f^{(m)}(x)|^2\dup x
  \end{equation}
  for all $f\in C_c^\infty(\Omega)$.
\end{thm}

\begin{proof}
  We prove the statement only for open intervals $(a,b)\subseteq\real$.
  Suppose that the above statement is true for some $m$. Then
  applying~\cite[Lemma 5.3.1]{Davi:1995} with $\alpha=-2m$,
  \begin{align*}
    \lefteqn{\frac{(2m+1)^2}4\int_a^b\frac{|f(x)|^2}{d(x)^{2(m+1)}}\dup x} 
    \qquad\qquad & \\
    &=\frac{(1+2m)^2}4\int_0^{(b-a)/2}x^{-2m-2}|f(x+a)|^2\dup x \\
    &\qquad\qquad\qquad\qquad+\frac{(1+2m)^2}4\int_0^
     {(b-a)/2}x^{-2m-2}|f(b-x)|^2\dup x \\
    &\le\int_0^{(b-a)/2}x^{-2m}|f'(x+a)|^2\dup x
     +\int_0^{(b-a)/2}x^{-2m}|f'(b-a)|^2\dup x \\
    &=\int_a^b\frac{|f'(x)|^2}{d(x)^{2m}}\dup x \\
    &\le\frac{4^m}{(2m-1)^2(2m-3)^2\dots1^2}\int_a^b|f^{(m+1)}(x)|^2\dup x.
  \end{align*}
  The first step of induction is dealt with 
  by~\cite[Corollary 5.3.2]{Davi:1995}.\finbox
\end{proof}

\begin{thm}{Lemma}\label{thm:groovy}
  Let $\omega\in S^{N-1}$. Then
  \begin{equation}
    \int_{S^{N-1}}\ip\xi\omega^{2m}\dup^{N-1}\omega=\frac{(2m-1)(2m-3)\dots1}
    {(N+2m-2)(N+2m-4)\dots N}|\xi|^{2m}.
  \end{equation}%
\end{thm}%
\begin{proof}
  Since the above integral is rotationally invariant and homogeneous of degree 
  $2m$ with respect to $\xi$ we see that
  \begin{equation*}
    \int_{S^{N-1}}\ip\xi\omega^{2m}\dup^{N-1}\omega=c|\xi|^{2m}.
  \end{equation*}
  Setting $\xi=(1,0,\dots,0)$ we see that for $N\ge3$
  \begin{align*}
    c&=\int_{S^{N-1}}\ip\xi\omega^{2m}\dup^{N-1}\omega \\
     &=\frac{1}{\omega_{N-1}}\int_{-\pi}^\pi\int_0^\pi\dots\int_0^\pi\cos^{2m}
       \theta_1\sin^{N-2}\theta_1\dots\sin\theta_{N-2}\dup\theta_1\dots\dup
       \theta_{N-1} \\
     &=\frac{(2m-1)(2m-3)\dots1}{(N+2m-2)(N+2m-4)\dots N},
  \end{align*}
  where $\omega_{N-1}$ denotes the surface area of the unit spherical shell 
  $S^{N-1}$ regarded as a subset of $\real^N$. The last step of this 
  calculation 
  requires elementary analysis and is therefore omitted. The cases $N=1,2$ 
  are simple.\finbox
\end{proof}

We may now prove the Hardy-Rellich operator inequality:

\begin{thm}{Theorem}\label{thm:hardy}
  Let $\poly=\power$ be the polyharmonic operator of order $2m$ acting in $L^2
  (\Omega)$, where $\Omega$ is a region in $\real^N$, and let $a_m$ be the 
  corresponding pseudodistance. Then, in the quadratic form sense, 
  \begin{equation}\label{eqn:myhardy}
    \poly\ge\frac{(N+2m-2)(N+2m-4)\dots N(2m-1)(2m-3)\dots1}{4^ma_m^{2m}}.
  \end{equation}
\end{thm}

\begin{proof}
  Let $f\in C_c^\infty(\Omega)$. Let $\omega\in S^{N-1}$ be fixed, and let $\{
  u_1=\omega,u_2,\dots,u_N\}$ be an orthonormal basis of $\real^N$. Let 
  $v=(v_1,\dots,v_N)$ denote coordinates with respect to that basis and let 
  $P$ be the coordinate transition matrix $x=vP$ from $v$ coordinates to 
  standard coordinates. Let $\hat v=(v_2,\dots,v_N)$ be fixed, and let $\Omega_
  {\hat{v}}$ be the open (not necessarily connected) set
  \begin{equation*}
    \Omega_{\hat{v}}=\{v_1\in\real:vP\in\Omega\}.
  \end{equation*}
  Define $g_{\hat{v}}$ and $d_{\hat{v}}$ by
  \begin{gather*}
  g_{\hat{v}}(v_1):=f(vP) \\
  d_{\hat{v}}(v_1):=d_\omega(vP).
  \end{gather*}
  Then $g_{\hat{v}}\in C_c^\infty(\Omega_{\hat{v}})$ and
  \begin{equation*}
    d_{\hat{v}}(v_1)=\min\{|y-v_1|:y\not\in\Omega_{\hat{v}}\}.
  \end{equation*}
  Using Lemma~\ref{thm:onedimhardy},
  \begin{equation*}
    \frac{(2m-1)^2(2m-3)^2\dots1^2}{4^m}\int_{\Omega_{\hat{v}}}\frac{|g_{\hat 
    v}(v_1)|^2}{d_{\hat v}(v_1)^{2m}}\dup v_1\le\int_{\Omega_{\hat{v}}}|g_{
    \hat v}^{(m)}(v_1)|^2\dup v_1,
  \end{equation*}
  and hence
  \begin{align*}
    \lefteqn{\frac{(2m-1)^2(2m-3)^2\dots1^2}{4^m}\int_\Omega\frac{|f(x)|^2}{a_m
    (x)^{2m}}\dup^Nx} \qquad & \\
    &=\frac{(2m-1)^2(2m-3)^2\dots1^2}{4^m}\int_{S^{N-1}}\int_\Omega\frac{|f(x)|
     ^2}{d_\omega(x)^{2m}}\dup^Nx\dup^{N-1}\omega \\
    &=\frac{(2m-1)^2(2m-3)^2\dots1^2}{4^m}\int_{S^{N-1}}\int_{\real^{N-1}}\int_
     {\Omega_{\hat{v}}}\frac{|g_{\hat v}(v_1)|^2}{d_{\hat v}(v_1)^{2m}}\dup v_1
     \dup^{N-1}\hat v\dup^{N-1}\omega \\
    &\le\int_{S^{N-1}}\int_{\real^{N-1}}\int_{\Omega_{\hat{v}}}|g_{\hat v}^{(m)
    }(v_1)|^2\dup v_1\dup^{N-1}\hat v\dup^{N-1}\omega \\
    &=\int_{S^{N-1}}\int_{\Omega}|\partial_\omega^mf(x)|^2\dup^Nx\dup^{N-1}
    \omega \\
    &=\int_{S^{N-1}}\int_{\real^N}\ip\xi\omega^{2m}|\hat f(\xi)|^2\dup^N
     \xi\dup^{N-1}\omega \\
    &=\frac{(2m-1)(2m-3)\dots1}{(N+2m-2)(N+2m-4)\dots N}\int_{\real^N}|\xi|^{2
    m}|\hat f(\xi)|^2\dup^N\xi \\
    &=\frac{(2m-1)(2m-3)\dots1}{(N+2m-2)(N+2m-4)\dots N}Q_m(f). \lowearlybox
  \end{align*}%
\end{proof}%
\begin{thm}{Corollary}\label{thm:conv}
  Suppose that $\Omega$ is a convex region in $\real^N$. Then
  \begin{equation}\label{eqn:convex}
    \poly\ge\frac{(2m-1)^2(2m-3)^2\dots1^2}{4^md^{2m}}.
  \end{equation}%
\end{thm}%
\begin{proof}
  Let $x\in\Omega$ and let $y\in\partial\Omega$ be such that $|y-x|=d(x)$.
  Suppose $z\in\Omega$. Constructing the point
  \begin{equation}
    p=y+\frac{\ip{z-y}{x-y}}{\ip{z-y}{z-y}}(z-y),
  \end{equation}
  we see that
  \begin{equation*}
    |p-x|^2=|y-x|^2-\frac{\ip{z-y}{x-y}\text{\makebox[0cm][l]{$^2$}}}{\ip{z-y}
    {z-y}}
  \end{equation*}
  so either
  \[ \ip{z-y}{x-y}=0 \]
  or
  \[ |p-x|<d(x). \]
  In the second case, $p$ will lie in $\Omega$, and so by convexivity the line
  segment joining $z$ and $p$ lies in $\Omega$. See Figure~\ref{fig:constr_p}.
  Since $y\not\in\Omega$, it cannot lie on this segment so
  \begin{equation*}
    \ip{z-y}{p-y}>0.
  \end{equation*}
  From the definition of $p$, this implies that
  \begin{equation*}
    \ip{z-y}{x-y}>0.
  \end{equation*}%
  \dbledgram{%
    \renewcommand{\afigwidth}{25}%
    \renewcommand{\afigheight}{31}%
    \renewcommand{\bfigwidth}{25}%
    \renewcommand{\bfigheight}{31}}{%
    \fig{hyper}{%
      \label{fig:constr_p}}{\put(78,103){$x$}%
      \put(188,93){$y$}%
      \put(51,229){$z$}%
      \put(142,150){$p$}%
      \put(9,17){$\Omega$}%
      \put(145,231){$H$}}{Construction of $p$}}{%
    \fig{hardy}{%
      \label{fig:reldist}}{\put(188,93){$y$}%
      \put(78,103){$x$}%
      \put(76,145){$d_\omega(x)$}%
      \put(109,80){$d(x)$}%
      \put(9,17){$\Omega$}%
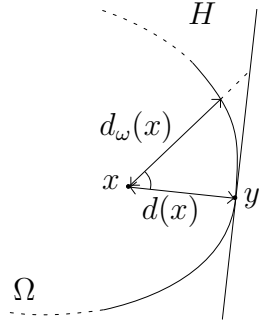
      \put(145,231){$H$}}{Relationships between distances}}%
  In both cases, $z$ lies in the set
  \[ \{z\in\real^N:\ip{z-y}{x-y}\ge0\}. \]
  Since $\Omega$ is open, it must therefore be a subset of the open half
  \begin{equation}
    H:=\{z\in\real^N:\ip{z-y}{x-y}>0\}
  \end{equation}
  of $\real^N$. From the definition~\eqref{eqn:defd} of $d_\omega$ we see 
  that
  \begin{align*}
    d_\omega(x)|\ip{y-x}\omega|&\le\min\{|s|:x+s\omega\not\in H\}\big|\bip{
    \frac{y-x}{|y-x|}}\omega\big|d(x) \\
    &=d(x)^2.
  \end{align*}
  See Figure~\ref{fig:reldist} for a diagrammatic representation of this last
  step. Hence
  \begin{equation*}
    |\ip{y-x}\omega|^{2m}d(x)^{-4m}\le d_\omega(x)
    ^{-2m}.
  \end{equation*}
  Therefore
  \begin{align}\label{eqn:convreg}
    \lefteqn{\frac{(2m-1)(2m-3)\dots1}{(N+2m-2)(N+2m-4)\dots N}d(x)^{-2m}} 
    \qquad\qquad\qquad\qquad & \notag\\
    &=\frac{(2m-1)(2m-3)\dots1}{(N+2m-2)(N+2m-4)\dots N}|y-x|^{2m}d(x)^{-4m} 
    \notag\\
    &=\int_{S^{N-1}}|\ip{y-x}\omega|^{2m}d(x)^{-4m}\dup^{N-1}\omega \notag\\
    &\le\int_{S^{N-1}}d_\omega(x)^{-2m}\dup^{N-1}\omega \notag\\
    &=a_m(x)^{-2m}.
  \end{align}
  Using Theorem~\ref{thm:hardy} we see that for $f\in C_c^\infty(\Omega)$,
  \begin{align*}
    \lefteqn{\frac{(2m-1)^2(2m-3)^2\dots1^2}{4^m}\int_\Omega\frac{|f(x)|^2}{d(x
    )^{2m}}\dup^Nx} \qquad & \\
    &\le\frac{(N+2m-2)(N+2m-4)\dots N(2m-1)(2m-3)\dots1}{4^m}\int_\Omega\frac{
    |f(x)|^2}{a_m(x)^{2m}}\dup^Nx \\
    &\le Q_m(f).\tag*{\finbox}
  \end{align*}%
\end{proof}%
\begin{note}
  It is simple to deduce a crude lower bound
  \[ \lambda_1\ge\frac{(2m-1)^2(2m-3)^2\dots1^2}{4^m\inradius(\Omega)^{2m}} \]
  on the first eigenvalue of $\poly$ for regions with finite inradius 
  \begin{equation}
    \inradius(\Omega):=\sup_{x\in\Omega}d(x).
  \end{equation}
  Since the strength of the inequality~\eqref{eqn:convex} lies in the values 
  of the potential near the boundary, the constant in the above bound is not 
  sharp.\finbox
\end{note}

\begin{note}
  The constants in Theorem~\ref{thm:hardy} and Corollary~\ref{thm:conv} are 
  optimal. This can be seen by choosing $\Omega=\{(x_1,\dots,x_N):x_1>0\}$ 
  and by considering the sequence of functions $f_n\in C_c^\infty(\Omega)$ 
  defined by
  \[ f_n(x)=x_1^{m-1/2}\phi_n(x_1)\psi(\hat x)  \]
  where $\psi\in C_c^\infty(\real^{N-1})$ and $\phi_n\in C_c^\infty((0,\infty))
  $ is chosen so that $\phi_n=1$ on the interval $[2/n,1]$, $\phi_n=0$ on $
  \real\setminus[1/n,2]$, $|D^j\phi_n|\le cn^j$ on $[1/n,2/n]$ and $|D^j\phi_n
  |\le c$ on $[1,2]$, for $j=0,1,\dots,2m$.

  Calculations now show that
  \[ \int_\Omega\frac{|f_n(x)|^2}{d(x)^{2m}}\dup^Nx\ge \int_{[2/n,1]\times\real
     ^{N-1}}x_1^{-1}|\psi(\hat{x})|^2\dup^Nx=\|\psi\|_2^2\ln n/2 \]
  and
  \[ Q_m(f_n)\le\frac{(2m-1)^2(2m-3)^2\dots1}{4^m}\|\psi\|_2^2\ln n/2+c'.  \]
  The constant in Corollary~\ref{thm:conv} is therefore optimal, and so the 
  constant in Theorem~\ref{thm:hardy} must also be optimal.\finbox
\end{note}
\section{Spectral Implications of the Inequality}%
\label{sec:impls}%
In the course of proving Corollary~\ref{thm:conv} we show, in 
inequality~\eqref{eqn:convreg}, that the pseudo\-distance $a_m$ is uniformly 
comparable to the boundary distance function $d$. This motivates the 
introduction of the following terminology:

\begin{definition}
  A region $\Omega$ is said to be regular if there is a constant $k<\infty$ 
  such that 
  \begin{equation}
    d(x)\le a_1(x)\le kd(x).
  \end{equation}
  for all $x\in\Omega$. More generally, we shall say that a region $\Omega$ is 
  $m$-regular if there is a constant $k_m<\infty$ such that
  \begin{equation}\label{eqn:mreg}
    d(x)\le a_m(x)\le k_md(x).
  \end{equation}
  for all $x\in\Omega$.
\end{definition}

\begin{definition}
  We shall say that $\Omega$ satisfies a uniform external ball condition if 
  there exist positive constants $\alpha$, $\beta$ such that for any $y\in
  \Omega$ and $0<s\le\beta$ there exists a ball $B(a;r)$ with center $a$ 
  satisfying $|a-y|\le s$, and radius $r$ satisfying $r\ge\alpha s$, which 
  does not meet $\Omega$.
\end{definition}

\begin{thm}{Examples}
  If any one of the following geometrical conditions is satisfied then the 
  region $\Omega\subseteq\real^N$ is regular:
  \begin{enumerate}
    \item $\Omega$ satisfies a uniform external ball condition with 
          $\beta=\infty$.
    \item $\Omega$ has finite inradius and satisfies a uniform external ball 
          condition.
    \item There exists a positive constant $c$ such that 
          \begin{equation*}
            |\{y\not\in\Omega:|y-a|<r\}|\ge cr^N
          \end{equation*}
          for all $a\in\partial\Omega$ and all $r>0$. 
  \end{enumerate}%
\end{thm}%
\begin{proof}
  See~\cite[Theorems 1.5.4 and 1.5.5]{Davi:1989} 
  and~\cite[Theorem 5.3.6]{Davi:1995}. The common characteristic of these 
  situations is that at any point $x\in\Omega$, the directional distance $d_\omega(x
  )$ to the boundary is uniformly comparable to the actual distance $d(x)$ to 
  the boundary over a uniform solid angle.\finbox
\end{proof}

\begin{thm}{Lemma}\label{thm:regular}
  Let $x\in\Omega$ be fixed. Then $a_m(x)$ is a decreasing function of $m$. 
  Hence if $\Omega$ is regular then it is $m$-regular for all $m\ge1$.
\end{thm}

\begin{proof}
  Let $\|.\|_p$ be norms on the spaces $L^p(S^{N-1},\dup^{N-1}\omega)$. Since the surface
  measure $\dup^{N-1}\omega$ in Definition~\ref{def:pseud} is normalized, 
  H\"older's inequality implies that for $m\ge n$
  \begin{equation*}
    \|d_\omega(x)^{-1}\|_{2n}\le\|d_\omega(x)^{-1}\|_{2m}\|1\|_{2mn/(m-n)}=\|d_\omega(x)^{-
    1}\|_{2m}.
  \end{equation*}
  Hence from Definition~\ref{def:pseud},
  \begin{equation*}
    a_m(x)=\|d_\omega(x)^{-1}\|_{2m}^{-1}\le\|d_\omega(x)^{-1}\|_{2n}^{-1}=a_n(x).
    \earlybox
  \end{equation*}
\end{proof}%
\begin{thm}{Theorem}
  Suppose that $\Omega$ is $m$-regular. Then $0\not\in\spec(\power
  )$ if and only if the inradius of $\Omega$ is finite.
\end{thm}

\begin{proof}
  Suppose that $\Omega$ has finite inradius. The Hardy-Rellich operator
  inequality~\eqref{eqn:myhardy} and $m$-regularity~\eqref{eqn:mreg} imply 
  that 
  \begin{equation}
    \power\ge\frac{(N+2m-2)(N+2m-4)\dots N(2m-1)(2m-3)\dots1}{4^mk_m^{2m}
    \inradius(\Omega)^{2m}}.
  \end{equation}
  Conversely, suppose that $d(x)$ is unbounded. For any $r>0$ there exists a
  ball $B_r$ with radius $r$, contained in $\Omega$. Using the Rayleigh-Ritz 
  variational formula (see~\cite[Section 4.5]{Davi:1995}),
  \begin{equation*}
    0\le\min(\spec(\power))\le\min(\spec(\power[B_r]))=r^{-2m}\min(\spec(
    \power[B_1])).
  \end{equation*}
  Hence $0\in\spec(\power)$.\finbox
\end{proof}

Note that to prove the above theorem one only needs the Hardy inequality for 
$m=1$ and inequality~\eqref{eqn:polylarger}. This approach, however, is not 
valid for a proof of the following theorem.

\begin{thm}{Theorem}\label{thm:compres}
  Suppose that $\Omega$ is $m$-regular. Then $\power^{-1}$ is compact 
  if and only if $d(x)\rightarrow0$ as $x\rightarrow\infty$.
\end{thm}

\begin{proof}
  Using the Hardy-Rellich inequality~\eqref{eqn:myhardy},
  \begin{align*}
    \power&\ge\frac12\power+\frac{(N+2m-2)(N+2m-4)\dots N(2m-1)(2m-3)\dots1}
    {2.4^ma_m^{2m}} \\
    &\ge\frac12\power[\real^N]+\frac{(N+2m-2)(N+2m-4)\dots N(2m-1)(2m-3)\dots
    1}{2^{2m+1}k_m^{2m}d^{2m}}
  \end{align*}
  as quadratic forms in $L^2(\real^N)$. The last operator in the above 
  inequality has compact resolvent because it is a Schr\"odinger operator 
  whose potential
  \[ V=\frac{(N+2m-2)(N+2m-4)\dots N(2m-1)(2m-3)\dots1}{2^{2m+1}k_m^{2m}d^{2m}
     } \]
  satisfies $V(x)\rightarrow\infty$ as $|x|\rightarrow\infty$. See 
  \cite[Theorem 12.5.5]{Mazy:1985}, although this result is proved with the 
  unnecessary restriction that $N<2m$. Simple modification of the proof 
  of~\cite[Theorem XIII.67]{Reed+Simo4:1978} yields the result without any 
  such restriction. It now follows that $\power^{-1}$ is compact. 

  Conversely suppose $d(x)$ does not converge to zero as $x\rightarrow\infty$,~
  $x\in\Omega$. Then there exist $r>0$ and a sequence of balls 
  $B_i\subseteq\Omega$, each with radius $r$. Let $\phi_i$ be the groundstate 
  of the operator $\power[B_i]$. Then
  \begin{align*}
    \ip{\phi_i}{\phi_j}&=\delta_{ij} \\
    \ip{\power\phi_i}{\phi_j}&=c\delta_{ij}
  \end{align*}
  where $c$ is independent of $i$, $j$. Using the Rayleigh-Ritz 
  formula of section~\cite{Davi:1995} we see that $\power^{-1}$ cannot be 
  compact.\finbox
\end{proof}%
\section{Lower Bound on the Trace of the Polyharmonic Semigroup}
\label{sec:lwbd}
In the remaining sections we build upon the methods of Davies~\cite{Davi:1985}
to obtain 
lower and upper bounds on the trace of the semigroup $e^{-\power t}$ and the 
resolvent $\power^{-\gamma}$.
The proof of the lower bound in Theorem~\ref{thm:lwtc} requires the following
sequence of lemmas:

\begin{thm}{Lemma}\label{thm:powereval}
  Let $\lambda_{m,n}$ denote the $n$-th eigenvalue of the polyharmonic operator
  $\poly$ acting in $L^2((0,1))$. Then
  \begin{equation}
    [n\pi]^{2m}\le\lambda_{m,n}\le[(m+n-1)\pi]^{2m}.
  \end{equation}%
\end{thm}%
\begin{proof}
  The left hand inequality is a consequence of 
  inequality~\eqref{eqn:polylarger}. We prove the other inequality as follows:
  Let $f_r\in W_0^{m,2}((0,1))$ be defined by
  \[ f_r(x)=\sin^{m-1}\pi x\sin r\pi x. \]
  Then $f_r\in M_{r+m-1}$ where
  \[ M_s=\lin\{1,\sin\pi x,\cos\pi x,\dots,\sin s\pi x,\cos s\pi x\}. \]
  Let $L_n\subseteq W_0^{m,2}([0,1])$ be defined by
  \[ L_n=\lin\{f_r:1\le r\le n\}. \]
  Then $L_n\subseteq M_{n+m-1}$, and by the Rayleigh-Ritz 
  formula~\cite{Davi:1995},
  \begin{align}\label{eqn:noth}
    \lambda_{m,n}&\le\sup\{Q_m(f):f\in L_n,\|f\|_2=1\} \notag\\
    &=\sup\left\{\|D^mf\|_2^2:f\in L_n,\|f\|_2=1\right\} \notag\\
    &\le\sup\left\{\|D^mf\|_2^2:f
    \in M_{n+m-1},\|f\|_2=1\right\}.
  \end{align}
  Suppose that $f\in M_s$ and $\|f\|_2=1$. Then
  \[ f(x)=\alpha_0+\sum_{r=1}^s(\alpha_r\sqrt2\cos r\pi x+\beta_r
  \sqrt2\sin r\pi x), \]
  where
  \[ \sum_{r=0}^s(\alpha_r^2+\beta_r^2)=1. \]
  Now
  \begin{equation*}
    \|D^mf\|_2^2=\sum_{r=1}^s(\alpha_r^2.(r\pi)^{2m}+\beta_r^2.(r\pi)^{2m})\le
    (s\pi)^{2m}.
  \end{equation*}
  Hence by inequality~\eqref{eqn:noth},
  \begin{equation*}
    \lambda_{m,n}\le[(m+n-1)\pi]^{2m}. \earlybox
  \end{equation*}%
\end{proof}%
\begin{thm}{Lemma}\label{thm:nomixed}
  The operator
  \begin{equation}\label{eqn:strop}
    H'=\sum_{i=1}^N\left(-\frac{\partial^2\phantom{x}}{\partial x_i^2}\right)^m
  \end{equation}
  acting in $L^2(C)$ with Dirichlet boundary conditions, where $C=(0,\delta
  )^N$, is uniformly elliptic, homogeneous of order $2m$, and has compact 
  resolvent. The eigenvalues of $H'$ are given by
  \begin{equation}
    \mu_n=\delta^{-2m}\sum_{i=1}^N\lambda_{m,n_i},
  \end{equation}
  where $n=(n_1,\dots,n_N)$ is a non-negative multi-index and $\lambda_{m,n_i}$
  are the eigenvalues of the one-dimensional polyharmonic operator in 
  Lemma~\slref{thm:powereval}.
\end{thm}

\begin{proof}
  Using the Fourier transform we may write the quadratic form $Q'$ of the 
  operator $H'$ as
  \begin{equation}
    Q'(f)=\int_{\real^N}\sum_{i=1}^N\xi_i^{2m}|\hat{f}(\xi)|^2\dup^N\xi.
  \end{equation}
  The symbol of $H'$ is
  \begin{equation}
    a(x,\xi)=\sum_{i=1}^N\xi_i^{2m},
  \end{equation}
  and is homogeneous of degree $2m$. Since
  \begin{equation}\label{eqn:itshomog}
    N^{-(m-1)}|\xi|^{2m}\le a(x,\xi)\le|\xi|^{2m},
  \end{equation}
  we see that $H'$ is uniformly elliptic.

  Let $\{f_{m,n}\}_{n=1}^\infty$ be the orthonormal sequence of eigenfunctions 
  corresponding to the eigenvalues $\lambda_{m,n}$ of $\poly$ acting in $L^2((0
  ,1))$. For each non-negative multi-index $n$ define
  \begin{equation}
    f_n(x)=\delta^{-N/2}\prod_{i=1}^Nf_{m,n_i}(x_i/\delta).
  \end{equation}
  Since the functions $f_{m,n}$ form a complete orthonormal set in $L^2((0,1))
  $, the functions $f_n$ form a complete orthonormal set in $L^2(C)$. 
  Moreover, since the $f_n$ are eigenfunctions of $H'$ with corresponding 
  eigenvalues $\mu_n$, we have found a complete list of eigenvalues.\finbox
\end{proof}

To find the lower bound on $\tc[e^{-\power t}]$ we use the technique of 
decomposing $\Omega$ into dyadic cubes by introducing Dirichlet boundary 
conditions along various internal partitioning surfaces. The cubes 
$C\subseteq\Omega$ we use are of the form
\begin{equation}\label{eqn:cubes}
  C=\left\{x\in\real^N:\frac{a_i}{2^n}<x_i<\frac{a_i+1}{2^n}\right\}
\end{equation}
for some $n\in\integer$ and some $a\in\integer^N$. 
Ordering the dyadic cubes~\eqref{eqn:cubes}
by inclusion, let $\{C_r:r\in\natural\}$ be an enumeration of the maximal 
cubes contained in $\Omega$, provided at least one exists. Let $\delta_r$ be 
the side length of $C_r$ and let $\Omega'=\bigcup_{r=1}^\infty C_r$.

\begin{thm}{Lemma}\label{thm:partition}
  The cubes $C_r$ are disjoint. Suppose that the inradius of $\Omega$ is 
  finite. Then $\overline{\Omega'}=\overline\Omega$, and moreover for $x\in
  \overline{C_r}$ we have $d(x)\le2N^{1/2}\delta_r$.
\end{thm}

\begin{proof}
  The inclusion $\overline{\Omega'}\subseteq\overline\Omega$ is obvious. 
  Conversely, suppose that $x\in\Omega$. Then $\overline{B(x;d(x))}\subseteq
  \overline\Omega$ and so $x$ will lie in some closed dyadic cube $\overline C
  \subseteq\overline\Omega$ with diameter at least $d(x)/2$. The edge length of 
  such a cube will be at least $d(x)/(2N^{1/2})$. Since $d(x)$ is bounded, the 
  point $x$ will lie in a maximal cube $\overline{C_r}$ with edge length 
  $\delta_r\ge d(x)/(2N^{1/2})$. Hence $\Omega\subseteq\overline{\Omega'}$ 
  and $d(x)\le2N^{1/2}\delta_r$.\finbox
\end{proof}

We shall always assume that $\Omega$ has finite inradius, for otherwise using 
Theorem~\ref{thm:compres}, we may deduce that $\tc[e^{-\power t}]=\infty$.

\begin{thm}{Theorem}\label{thm:lwtc}
  For $0<t<\infty$
  \begin{equation}\label{eqn:lwtc}
    b_{m,N}t^{-N/2m}\int_\Omega \exp[{-c_{m,N}d^{-2m}t}]\le\tc[e^{-\power t}],
  \end{equation}
  where 
  \begin{equation*}
    b_{m,N}=N^{-N(m-1)/2m}(2\pi)^{-N}\Gamma(1+1/2m)2^{N/2m}
  \end{equation*}
  and
  \begin{equation*}
    c_{m,N}=(4mN\pi)^{2m}/2.
  \end{equation*}%
\end{thm}%
\begin{proof}
  Let $\Omega=(0,1)\subseteq\real$. Using the notation of 
  Lemma~\ref{thm:powereval} and the spectral mapping theorem, we see that
  the trace of the semigroup $e^{-\power[{[0,1]}]t}$ is $\sum_{n=1}^\infty e^{-
  \lambda_{m,n}t}$, and moreover that
  \begin{align}\label{eqn:traoned}
    \sum_{n=1}^\infty e^{-\lambda_{m,n}t}&\ge\sum_{n=0}^\infty e^{-[(m+n)\pi]^
    {2m}t} \notag\\
    &\ge\sum_{n=0}^\infty e^{-2^{2m-1}(m^{2m}+n^{2m})\pi^{2m}t} \notag\\
    &\ge e^{-(2m\pi)^{2m}t/2}\int_0^\infty e^{(2x\pi)^{2m}t/2}\dup x \notag\\
    &=[(2\pi)^{2m}t/2]^{-1/2m}\Gamma(1+1/2m)e^{-(2m\pi)^{2m}t/2} \notag\\
    &=b_{m,1}t^{-1/2m}\exp[-c_{m,1}2^{-2m}t].
  \end{align}
  Let $H'$ denote the operator~\eqref{eqn:strop} acting in $L^2(C_r)$ with 
  Dirichlet boundary conditions. Inequality~\eqref{eqn:itshomog} implies that 
  \begin{equation*}
    \power[C_r]\le N^{m-1}H',
  \end{equation*}
  and hence using Lemma~\ref{thm:nomixed} and 
  equation~\eqref{eqn:traoned},
  \begin{align*}
    \tc[e^{-\power[C_r]t}]&\ge\tc[e^{-N^{m-1}H't}] \\
    &=\sum_{n\in\natural^N}e^{-N^{m-1}\delta_r^{-2m}t\sum_{i=1}^N\lambda_
    {m,n_i}} \\
    &=\sum_{n\in\natural^N}\prod_{i=1}^Ne^{-N^{m-1}\delta_r^{-2m}t\lambda_
    {m,n_i}} \\
    &=\bigg[\sum_{n=1}^\infty e^{-N^{m-1}\delta_r^{-2m}t\lambda_{m,n}}\bigg]
    ^N \\
    &\ge b_{m,1}^N[N^{m-1}\delta_r^{-2m}t]^{-N/2m}\exp[-c_{m,1}2^{-2m}N^m
    \delta_r^{-2m}t] \\
    &=b_{m,N}\delta_r^Nt^{-N/2m}\exp[-c_{m,N}(2N^{1/2}\delta_r)^{-2m}t].
  \end{align*}
  Using the results of Lemma~\ref{thm:partition},
  \begin{align*}
    \tc[e^{-\power t}]&\ge\tc[e^{-\power[\Omega']t}] \\
    &=\sum_{r=1}^\infty\tc[e^{-\power[C_r]t}] \\
    &\ge\sum_{r=1}^\infty b_{m,N}\delta_r^Nt^{-N/2m}\exp[-c_{m,N}(2N^{1/2}
    \delta_r)^{-2m}t] \\
    &=b_{m,N}t^{-N/2m}\sum_{r=1}^\infty\int_{C_n}\exp[-c_{m,N}td(x)^{-2m}]\dup^Nx 
    \\
    &=b_{m,N}t^{-N/2m}\int_\Omega\exp[-c_{m,N}td^{-2m}]. \lowearlybox
  \end{align*}
\end{proof}%
\section{Upper Bound on the Trace of the Polyharmonic Semigroup}
\label{sec:upbd}
In order to prove an upper bound on the trace we shall need to assume that the 
region $\Omega$ satisfies the following condition.

\begin{thm}{Condition}\label{thm:bdker}
  Let $\Omega$ be a region such that the kernel $K_\Omega(t,x,y)$ of 
  $e^{-\power t}$ exists, is jointly continuous and satisfies
  \begin{equation}
    |K_\Omega(t,x,y)|\le ct^{-N/2m}. 
  \end{equation}
  for some $c=c_\Omega$, and for all $t>0$ and $x,y\in\Omega$. Let $b_{m,N}'=2^
  {N/2m}c$.
\end{thm}

Two special cases in which this condition is satisfied are given in the 
following two examples:

\begin{thm}{Example}\label{thm:lapkerbd}
  For all $N$ the Laplacian $(-\Delta)\dir$ acting in $L^2(\Omega)$ has a heat 
  kernel $K(t,x,y)$ which satisfies
  \[ 0\le K(t,x,y)\le(4\pi t)^{-N/2} \]
  for all $x,y\in\Omega$ and $t>0$.
\end{thm}

\begin{proof}
  See~\cite[example 2.1.8]{Davi:1989}.\finbox
\end{proof}

\begin{thm}{Example}\label{thm:polykerbd}
  Suppose that $\Omega\subseteq\real^N$ and $N<2m$. Then $\poly$ acting in $L^2
  (\Omega)$ has a heat kernel which satisfies
  \begin{equation*}
    |K(t,x,y)|\le ct^{-N/2m}
  \end{equation*}
  for all $x,y\in\Omega$ and all $t>0$.
\end{thm}

\begin{proof}
  By the spectral mapping theorem we see that
  \[ \|\power^{1/2}e^{-\power t}\|_{2,2}\le ct^{-1/2}. \]
  For $f\in L^2(\Omega)$ and $t>0$, let $f_t=e^{-\power t}f\in W_0^{m,2}(
  \Omega)$. Using a standard Sobolev embedding theorem,
  \begin{align*}
    \|f_t\|_\infty&\le c\|\power^{1/2}f_t\|^{N/2m}\|f_t\|_2^{1-N/2m} \\
    &\le c\|\power^{1/2}e^{-\power t}\|_{2,2}^{N/2m}\|f\|_2^{N/2m}\|f\|_2^{1-N/
    2m} \\
    &\le ct^{-N/4m}\|f\|_2.
  \end{align*}
  Hence
  \begin{equation*}
    \|e^{-\power t}\|_{\infty,2}\le ct^{-N/4m}.
  \end{equation*}
  By duality,
  \begin{equation}\label{eqn:l1libd}
    \|e^{-\power t}\|_{\infty,1}\le\|e^{-\power t}\|_{\infty,2}\|e^{-\power t}
    \|_{2,1}\le ct^{-N/2m}.
  \end{equation}
  Define $\phi_t:\Omega\rightarrow L^2(\Omega)\cap L^\infty(\Omega)$ by the 
  property
  \begin{equation*}
    \ip f{\phi_t(x)}=(e^{-\power t}f)(x)
  \end{equation*}
  for all $f\in L^1(\Omega)\cap L^2(\Omega)$. Since $e^{-\power t}f$ is a 
  smooth function, the map $x\mapsto\phi_t(x)$ from $\Omega\rightarrow L^2(
  \Omega)$ is smooth in the weak Hilbert space sense and hence, 
  by~\cite[Corollary 1.42]{Davi:1980}, it is smooth. Define
  \begin{equation*}
    K(t,x,y)=[\phi_t(x)](y).
  \end{equation*}
  Then by the definition of $\phi_t$, we see that $K(t,.,.)$ is an integral
  kernel of $e^{-\power t}$. Using the identity
  \begin{equation*}
    [\phi_{s+t}(x)](y)=\ip{\phi_s(x)}{\phi_t(y)}
  \end{equation*}
  for all $t,s>0$ we see that $K$ is smooth in $x$ and $y$. Moreover, 
  by~\eqref{eqn:l1libd} we see that
  \begin{equation*}
    |K(t,x,y)|\le ct^{-N/2m}. \earlybox
  \end{equation*}%
\end{proof}%
\begin{thm}{Theorem}\label{thm:uptc}
  Suppose that $\Omega$ satisfies Condition~\slref{thm:bdker}. Then
  \begin{equation}\label{eqn:uptc}
    \tc[e^{-\power t}]\le b_{m,N}'t^{-N/2m}\int_\Omega\exp[-c_{N,m}'a_m
    ^{-2m}t]
  \end{equation}
  where $b_{m,N}'$ is determined by Condition~\slref{thm:bdker}, and
  \begin{equation*}
    c_{m,N}'=2^{-2m-1}(N+2m-2)(N+2m-4)\dots N(2m-1)(2m-3)\dots1.
  \end{equation*}%
\end{thm}%
\begin{proof}
  The Hardy-Rellich inequality~\eqref{eqn:myhardy} shows that
  \[ \power\ge\frac12\power+\frac{(N+2m-2)(N+2m-4)\dots N(2m-1)
     (2m-3)\dots1}{2.4^ma_m^{2m}}. \]
  Using the Golden-Thompson inequality~\cite{Lena:1971}, 
  integration of the kernel along the 
  diagonal~\cite[pages 65,66]{Reed+Simo3:1979}, and Condition~\ref{thm:bdker} 
  we see that
  \begin{align*}
    \tc[e^{-\power t}]&\le\tc[\exp[{-\power t/2-c_{N,m}'a_m^{-2m}t}]] \\
    &\le\tc[e^{-\power t/4}\exp[{-c_{N,m}'a_m^{-2m}t}]e^{-\power t/4}] \\
    &=\int_\Omega K_\Omega(t/2,x,x)\exp[-c_{N,m}'a_m(x)^{-2m}t]\dup^Nx \\
    &\le b_{m,N}'t^{-N/2m}\int_\Omega\exp[-c_{N,m}'a_m(x)^{-2m}t]\dup^Nx. 
    \lowearlybox
  \end{align*}
\end{proof}%
\section{Equivalent Conditions for Finite Trace}

We can now use the lower and upper bounds of Theorems~\ref{thm:lwtc} 
and~\ref{thm:uptc} to give conditions for finite trace of $e^{-\power t}$ and 
$\power^{-\gamma}$ in terms of integrals involving the distance function $d$.

\begin{thm}{Theorem}\label{thm:equiv}
  Suppose that $\Omega$ is $m$-regular and satisfies 
  Condition~\slref{thm:bdker}. Then
  \[ \tc[e^{-\power t}]<\infty \]
  for all $t\in(0,\infty)$ if and only if
  \[ \int_\Omega e^{-td^{-2m}}<\infty \]
  for all $t\in(0,\infty)$.
\end{thm}

\begin{proof}
  Since $\Omega$ is $m$-regular, inequality~\eqref{eqn:mytrace} becomes
  \begin{equation*}
    bt^{-N/2m}\int_\Omega\exp[-ctd^{-2m}]\le\tc[e^{-\power
    t}]\le b't^{-N/{2m}}\int_\Omega\exp[-c'tk_md^{-2m}].
  \end{equation*}
  \finbox
\end{proof}

\begin{thm}{Corollary}
  Suppose that $\Omega$ is $m$-regular and satisfies 
  Condition~\slref{thm:bdker}, and that $\gamma>N/2m$. Then
  \begin{equation*}
    \tc[\power^{-\gamma}]<\infty
  \end{equation*}
  if and only if
  \begin{equation*}
    \int_\Omega d^{2m\gamma-N}<\infty.
  \end{equation*}%
\end{thm}%
\begin{proof}
  Using Fubini's theorem for traces we see that
  \begin{align*}
    \int_0^\infty\tc[e^{-\power t}]t^{\gamma-1}\dup t&=\tc\bigg[\int_0^\infty 
    e^{-\power t}t^{\gamma-1}\dup t\bigg] \\
    &=\Gamma(\gamma)\tc[\power^{-\gamma}]
  \end{align*}
  Integration of inequality~\eqref{eqn:uptc} gives
  \begin{align*}
    \int_0^\infty\tc[e^{-\power t}]t^{\gamma-1}\dup t&\le\int_\Omega\int_0^
    \infty b't^{-N/2m+\gamma-1}\exp[-c'ta_m(x)^{-2m}]\dup t\dup^Nx \\
    &=b'c'{}^{-\gamma+N/2m}\Gamma(\gamma-N/2m)\int_\Omega a_m^{2m
    \gamma-N}
  \end{align*}
  and similarly by integrating inequality~\eqref{eqn:lwtc} we see that
  \begin{equation*}
    bc^{-\gamma+N/2m}\Gamma(\gamma-N/2m)\int_\Omega d^{2m\gamma
    -N}\le\int_0^\infty\tc[e^{-\power t}]t^{\gamma-1}\dup t.
  \end{equation*}
  The result follows as in Theorem~\ref{thm:equiv} because $\Omega$ is 
  $m$-regular.\finbox
\end{proof}%

\bibliographystyle{plain}

\scshape
M. Owen: Erwin Schr\"odinger International Institute for Mathematical Physics,
Boltzmanngasse 9, A-1090 Wien, Austria.

\slshape
E-mail address: mowen@@esi.ac.at
\end{document}